\def\tpt{\hskip 20pt}

\input amstex
\input amsppt.sty
\magnification=\magstep1
\hsize=15truecm
\vsize=20truecm
\TagsOnRight
\pageno=1

{\tpt}

\vskip 20pt

\centerline {\bf A Geometric Proof of Mordell's Conjecture for
Function Fields\footnote""{Mathematics Subject Classification
(1991): 14H05, 14H10, 14G27}}

\vskip 10pt \centerline {Kezheng Li\footnote"*"{Supported by NSFC
of China, grant number 10531060}}

\vskip 10pt
\centerline{\sl Department of Mathematics, Capital
Normal University, Beijing 100037, China}


\vskip 20pt
\flushpar{\bf Abstract.} Let $\Cal C,\Cal C'$ be
curves over a base scheme $S$ with $g(\Cal C)\ge 2$. Then the
functor $T\mapsto\{$generically smooth $T$-morphisms
$T\times_S\Cal C'\to T\times_S\Cal C\}$ from $((S$-schemes)) to
((sets)) is represented by a quasi-finite unramified $S$-scheme.
From this one can deduce that for any two integers $g\ge 2$ and
$g'$, there is an integer $M(g,g')$ such that for any two curves
$C,C'$ over any field $k$ with $g(C)=g$, $g(C')=g'$, there are at
most $M(g,g')$ separable $k$-morphisms $C'\to C$. It is
conjectured that the arithmetic function $M(g,g')$ is bounded by a
linear function of $g'$.

\vskip 20pt
\flushpar {\bf 0. Introduction}

\vskip 10pt
We recall the works of Y. Manin, H. Grauert and P. Samuel on Mordell's conjecture for
function fields (see [G] and [S]). A main part of the conjecture can be stated as
follows.

\vskip 10pt
\flushpar{\bf Theorem 1}. Let $C,C'$ be two smooth projective curves over
a field $k$, where $C$ had genus $g\ge 2$. Then there are at most a finite
number of finite separable $k$-morphisms from $C'$ to $C$.

\vskip 10pt
From this one can deduce that

\vskip 10pt
\flushpar{\bf Theorem 2}. Let $C$ be a smooth projective curve over
a field $k$ with genus $g\ge 2$. Then for any finitely generated field
extension $K\supset k$, $C$ has at most a finite number of smooth
$K$-points over $k$.

\vskip 10pt
\flushpar Here a smooth $K$-point over $k$ means a smooth $k$-morphism
$\text{Spec}(K)\to C$, and this is equivalent to a $k$-algebra homomorphism
$\phi :K(C)\to K$ (where $K(C)$ is the function field of $C$) such that
$K\supset\text{im}(\phi )$ is a separably generated extension. Theorem 2
can be restated as

\vskip 10pt
\flushpar{\bf Theorem 3}. Let $C$ be a smooth projective curve over
a field $k$ with genus $g\ge 2$, and $X$ be a variety over $k$. Then there
are at most a finite number of generically smooth $k$-morphisms $X\to C$.

\vskip 10pt
\flushpar (See also [B], [Hir], [Hr], [N], [Vol], [Voj] for some recent
developments along this line.)

In this paper we will give a geometric proof for the above facts. Our main result
is

\vskip 10pt
\flushpar{\bf Theorem 4}. Let $S$ be a noetherian scheme and $\Cal C,\Cal C'$ be
curves over $S$ (i.e. smooth projective morphisms $\Cal C\to S$, $\Cal C'\to S$ of
relative dimension 1 with geometrically integral fibers). Suppose the fibers of
$\Cal C$ over $S$ all have genus $g\ge 2$. Then there is a quasi-finite unramified
$S$-scheme $\Cal T$ representing the following functor

\vskip 5pt
$((S$-schemes$))\to$((sets))

\hskip 53pt $T\mapsto\{$generically smooth $T$-morphisms
$f:T\times_S\Cal C'\to T\times_S\Cal C\}$

\vskip 5pt
\flushpar where ``generically smooth'' means there is an open subscheme
$U\subset T\times_S\Cal C'$, faithfully flat over $T$, such that $f|_U$ is smooth.
In particular, there are at most a finite number of generically smooth $S$-morphisms
from $\Cal C'$ to $\Cal C$.

\vskip 10pt
Some main ideas of the proof of Theorem 4 come from [L1] and [L2].

A special case of Theorem 4 is

\vskip 10pt
\flushpar{\bf Theorem 5}. Let $C,C'$ be smooth projective curves over a field
$k$ with genera $g(C)=g\ge 2$, $g(C')=g'$. Then there is a finite \'etale
$k$-scheme $\Cal T$ representing the following functor

\vskip 5pt
$((k$-schemes$))\to$((sets))

\hskip 52pt $T\mapsto\{$generically smooth $T$-morphisms
$f:T\times_kC'\to T\times_kC\}$

\vskip 10pt
This is slightly stronger than Theorem 1.

Using Theorem 4 and moduli theory on curves (cf. [M1] or [M2]), we can get

\vskip 10pt
\flushpar{\bf Theorem 6}. For any two positive integers $g\ge 2$ and $g'$,
there is an integer $M(g,g')$ such that for any two curves $C,C'$ over any field
$k$ with $g(C)=g$, $g(C')=g'$, it holds that

\vskip 5pt
$\#\{$finite separable $k$-morphisms $C'\to C\}\le M(g,g')$.

\vskip 5pt
\flushpar Furthermore, for any finitely generated field extension $K\supset k$, there is an
integer $M(g,K/k)$ such that for any curve $C$ over $k$ with genus $g$, it holds that

\vskip 5pt
$\#\{$smooth $K$-points of $C$ over $k\}\le M(g,K/k)$.

\vskip 5pt
\flushpar and for any variety $X$ over $k$, there is an
integer $M(g,X/k)$ such that for any curve $C$ over $k$ with genus $g$,

\vskip 5pt
$\#\{$generically smooth $k$-morphisms $X\to C\}\le M(g,X/k)$.

\vskip 10pt
This strengthens Theorem 1, Theorem 2 and Theorem 3.

From Theorem 4 and Theorem 6, one can also deduce that

\vskip 10pt
\flushpar{\bf Theorem 7}. Let $S$ be a noetherian scheme and $\Cal C_1,...,\Cal C_m$,
$\Cal C_1^{\prime},...,\Cal C_n^{\prime}$ be curves over $S$, such that the fibers of
$\Cal C_i$ over $S$ all have genus $g_i\ge 2$ ($1\le i\le m$), and the fibers of
$\Cal C_j^{\prime}$ over $S$ all have genus $g_j^{\prime}$ ($1\le j\le n$). Let
$X=\Cal C_1\times_S\cdots\times_S\Cal C_m$,
$Y=\Cal C_1^{\prime}\times_S\cdots\times_S\Cal C_n^{\prime}$.

\vskip 5pt\leftskip=30pt\parindent=-10pt
i) There is a quasi-finite unramified
$S$-scheme $\Cal T$ representing the following functor

\vskip 5pt
\flushpar $((S$-schemes$))\to$((sets))

\flushpar \hskip 53pt $T\mapsto\{$generically smooth $T$-morphisms
$f:T\times_SY\to T\times_SX\}$

\vskip 5pt
ii) If $S=\text{Spec}(k)$ for a field $k$, then any generically
smooth $k$-morphism
$f:Y\to X$ can be factored to a product of morphisms of curves. To be precise,
there is an injective map $\lambda :\{ 1,2,...,m\}\to\{ 1,2,...,n\}$ and $m$
finite separable morphisms $f_i:\Cal C_{\lambda (i)}^{\prime}\to\Cal C_i$ over $k$,
such that $f(x_1,...,x_n)=(f_1(x_{\lambda (1)}),...,f_m(x_{\lambda (m)}))$.

iii) Let $I$ be the set of all injective maps from $\{ 1,2,...,m\}$
to $\{ 1,2,...,n\}$. If $S=\text{Spec}(k)$ for a field $k$, then

\vskip 5pt
\flushpar $\#\{$generically smooth $k$-morphisms $Y\to X\}\le\sum\limits_{\lambda\in I}
\prod\limits_{i=1}^mM(g_i,g_{\lambda (i)}^{\prime})$.

\vskip 5pt
iv) Suppose $S=\text{Spec}(k)$ for a field $k$. Let $K=K(Y)$ (the function
field of $Y$). Then

\vskip 5pt
\flushpar $\#\{$smooth $K$-points of $X$ over $k\}\le\sum\limits_{\lambda\in I}
\prod\limits_{i=1}^mM(g_i,g_{\lambda (i)}^{\prime})$.

\vskip 10pt\leftskip=0pt\parindent=20pt
We conjecture that $M(g,g')$ is bounded by a linear function of $g'$. This holds
when $g=g'$.

\vskip 10pt
\flushpar{\bf Acknowledgement}. I wish to thank Kejian Xu for his stimulating
discussions with me on this topic. This paper is related to a recent work of his
(see [X]).

\vskip 30pt
\flushpar {\bf 1. The main theorem}

\vskip 10pt
We first fix some terminologies. Let $S$ be a scheme. By a {\sl curve} over $S$
we will mean a smooth projective morphisms $\pi:\Cal C\to S$ of relative
dimension 1 with geometrically integral fibers; in this case if all of the
fibers of $\pi$ have the same genus $g$, we say $\Cal C$ (or $\pi$) has {\sl genus} $g$.
Let $X,Y$ be schemes over $S$. An $S$-morphism $f:X\to Y$ is
called {\sl generically smooth} if there is an open subscheme $U\subset X$,
faithfully flat over $S$, such that $f|_U:U\to Y$ is smooth. If $K$ is a field
and $\text{Spec}(K)$ is an $S$-scheme, then an $S$-morphism $f:\text{Spec}(K)\to X$
is called a $K$-{\sl point} of $X$ over $S$; in this case we say the $K$-point
is {\sl smooth} if $f$ is smooth. If $X$ is a variety over a field $k$,
we denote by $K(X)$ the function field of $X$.

Note that when $X$ is a variety over a field $k$ and $K$ is a finitely generated
field extension
of $k$, a smooth $K$-point of $X$ over $k$ is equivalent to a $k$-algebra
homomorphism $\phi :K(X)\to K$ such that
$K\supset\text{im}(\phi )$ is a separably generated extension.

\vskip 10pt
\flushpar{\bf Lemma 1}. Let $S$ be a noetherian scheme and $X,Y$ be
flat projective $S$-schemes. Then there is a locally
quasi-projective $S$-scheme $\Cal T$ representing the following functor

\vskip 5pt
\hskip -10pt $\goth M\goth o\goth r_{\text{gs}}:((S$-schemes$))\to$((sets))

\hskip 78pt $T\mapsto\{$generically smooth $T$-morphisms
$T\times_SY\to T\times_SX\}$

\vskip 10pt
\flushpar {\sl Proof.} By moduli theory, there is a locally quasi-projective
$S$-scheme $\Cal T'$ representing the following functor

\vskip 5pt
$((S$-schemes$))\to$((sets))

\hskip 53pt $T\mapsto\{T$-morphisms $T\times_SY\to T\times_SX\}$

\vskip 5pt
\flushpar Here $\Cal T'$ is a locally closed subscheme of the Hilbert scheme
$\Cal Hilb_{Y\times_SX/S}$, which represents all of the flat closed subschemes
of $Y\times_SX$ over $S$. Let $f:\Cal T'\times_SY\to\Cal T'\times_SX$ be the
universal morphism over $\Cal T'$. Then there is a largest open subscheme
$U\subset\Cal T'\times_SY$ such that $f|_U$ is smooth. Let $\Cal T$
be the image of $U$ under the projection $\text{pr}_1:\Cal T'\times_SY\to\Cal T'$
(this makes sense because $\text{pr}_1$ is flat, hence is an open map).
It is easy to see that $\Cal T$ represents
$\goth M\goth o\goth r_{\text{gs}}$. \ \ \ Q.E.D.

\vskip 10pt
\flushpar{\bf Lemma 2}. Let $C,C'$ be two curves over a field $k$ with $g(C)=g\ge 2$,
$g(C')=g'$. Let $f:C'\to C$ be a separable morphism over $k$. Then $d=\deg (f)\le g'-1$.
Let $p'\in C'$, $p\in C$ be $k$-points. Let $X=C'\times_kC$, and take the ample
invertible sheaf of $X$ to be $O_X(D)$ for the divisor $D=p'\times_kC+C'\times_kp$.
Then the graph $\Gamma_f\subset X$ of $f$ has Hilbert polynomial
$\chi (x)=(d+1)x+1-g'$.

\vskip 10pt
\flushpar {\sl Proof.} By Hurwitz's theorem, we have
$2g'-2\ge d(2g-2)\ge 2d$, hence $d=\le g'-1$. Since $\Gamma_f\cong C'$,
by Riemann-Roch theorem we have $h^0(D')-h^1(D')=\deg (D')+1-g'$ for any
divisor $D'$ on $\Gamma_f$. Take $D'=D\cap\Gamma_f$, it is easy to
see that $\deg (D')=d+1$. Hence $\chi (n)=h^0(nD')-h^1(nD')=\deg (nD')+1-g'=(d+1)n+1-g'$.
This shows that $\chi (x)=(d+1)x+1-g'$. \ \ \ Q.E.D.

\vskip 10pt
\flushpar{\bf Proposition 1}. Let $S$ be a noetherian scheme and $\Cal C,\Cal C'$ be
curves over $S$ of genera $g\ge 2$ and $g'$ respectively. Let $\Cal T$ be the locally
quasi-projective $S$-scheme representing the following functor

\vskip 5pt
$((S$-schemes$))\to$((sets))

\hskip 53pt $T\mapsto\{$generically smooth $T$-morphisms
$f:T\times_S\Cal C'\to T\times_S\Cal C\}$

\vskip 10pt
\flushpar (as in Lemma 1). Then $\Cal T$ is quasi-finite and unramified over $S$.

\vskip 10pt
\flushpar {\sl Proof.}  First note that $\Cal Hilb_{\Cal C'\times_S\Cal C/S}$ is a
disjoint union of projective $S$-schemes $\Cal H^{\chi}$ indexed by
(infinitely many) Hilbert polynomials $\chi$. By Lemma 2, we see that
$\Cal T\subset\Cal Hilb_{\Cal C'\times_S\Cal C/S}$ is contained in a union
of a finite number of $\Cal H^{\chi}$s (for $\chi (x)=(d+1)x+1-g'$, $1\le d<g'$).
Hence $\Cal T$ is quasi-projective. Therefore it is enough to show that the
projection $q:\Cal T\to S$ is unramified.

Denote by $\Phi :\Cal T\times_S\Cal C'\to\Cal T\times_S\Cal C$ the universal morphism
over $\Cal T$, and $\rho =\text{pr}_2\circ\Phi :\Cal T\times_S\Cal C'\to\Cal C$.

Case 1: $S=\text{Spec}(k)$ for an algebraically closed field $k$.
Let $t\in\Cal T$ be a $k$-point, and let $f:\Cal C'\to\Cal C$ be the corresponding
separable $k$-morphism. Denote by $\zeta :\{ t\}\cong\text{Spec}(k)\to\Cal T$ the
inclusion. Let
$$\alpha =(\rho ,\rho\circ ((\zeta\circ\text{pr}_1)\times_k\text{id}_{C'}):
\Cal T\times_S\Cal C'\to\Cal C\times_k\Cal C \tag 1$$
i.e. $\alpha (t',x)=(\rho(t',x),f(x))$ ($\forall t'\in\Cal T$). We have a
commutative diagram
$$\CD \Cal C'\cong\  @. \text{Spec}(k)\times_k\Cal C' @>{f}>> \Cal C \cr
@. @VV{\zeta\times_k\text{id}_{\Cal C'}}V @VV{\Delta}V \cr @. \Cal
T\times_S\Cal C' @>{\alpha}>> \Cal C\times_k\Cal C \endCD\tag 2$$
Let $\Cal I$, $\Cal J$ and $\Cal J_0$ be the ideal sheaves of the
closed immersions $\Delta$, $\zeta\times_k\text{id}_{\Cal C'}$ and
$\zeta$ respectively. Clearly $\Cal J\cong\text{pr}_1^*\Cal J_0$.
By (2), $\alpha$ induces a morphism $\alpha^*\Cal I\to\Cal J$.
Applying $(\zeta\times_k\text{id}_{\Cal C'})^*$ we get a
homomorphism of coherent sheaves on $\Cal C'$:
$$\aligned \eta : & (\zeta\times_k\text{id}_{\Cal C'})^*(\alpha^*\Cal I)
\cong f^*(\Delta^*\Cal I)\cong f^*\Omega_{\Cal C/k}^1\to \cr
& (\zeta\times_k\text{id}_{\Cal C'})^*\Cal J\cong (\zeta^*\Cal J_0)\otimes_kO_{\Cal C'}
\cong O_{\Cal C'}^n \endaligned\tag 3$$
where $n=\dim_k(\zeta^*\Cal J_0)$.

We claim that $\eta =0$. Indeed, if $\eta\ne 0$, then there would
be a non-zero homomorphism $\eta':f^*\Omega_{\Cal C/k}^1\to
O_{\Cal C'}$. Since $f^*\Omega_{\Cal C/k}^1$ is an invertible
sheaf, $\eta'$ would be a monomorphism. Therefore we would have a
monomorphism $H^0(\Omega_{\Cal C/k}^1)\hookrightarrow
H^0(f^*\Omega_{\Cal C/k}^1)\hookrightarrow H^0(O_{\Cal C'})\cong
k$, contrary to $\dim_k(H^0(\Omega_{\Cal C/k}^1))=g\ge 2$.

Let $T\subset\Cal T$ be the closed subscheme defined by the ideal sheaf $\Cal J_0^2$,
and $V\subset C\times_k\Cal C$ be the closed subscheme defined by the ideal sheaf
$\Cal I^2$.
Then $\alpha$ induces a morphism $\alpha_1:T\to V$. By $\eta =0$ we see that
$\alpha_1^*(\Cal IO_V)=0$, hence $\alpha_1$ factors through $\Delta (\Cal C)$,
i.e. $\rho|_{T\times_k\Cal C'}=f\circ\text{pr}_2:T\times_k\Cal C'\to\Cal C$. By the
universality of $\Cal T$, we have a commutative diagram
$$\CD T @>{q}>> \text{Spec}(k) \cr
@VV{\text{inclusion}}V @VV{\zeta}V \cr
\Cal T @>{\text{id}}>> \Cal T \endCD\tag 4$$
where $q$ is the projection.
Hence $\Cal J_0O_T=0$, i.e, $\Cal J_0=\Cal J_0^2$. This means that
$$(\Omega_{\Cal T/k}^1)_t\cong\zeta^*\Cal J_0=0 \tag 5$$
where $(\Omega_{\Cal T/k}^1)_t$ is the fiber of $\Omega_{\Cal T/k}^1$ at $t$,
and $\zeta^*\Cal J_0$ can be viewed as $\Cal J_0/\Cal J_0^2$ restricted to $\{ t\}$.
Since $t$ is an arbitrary closed
point of $\Cal T$, we have $\Omega_{\Cal T/k}^1=0$,
i.e. $\Cal T$ is unramified over $k$.

Case 2: $S=\text{Spec}(k)$ for an arbitrary field $k$. Let $\bar k$ be the algebraic
closure of $k$, and denote $\bar{\Cal C}=\Cal C\otimes_k\bar k$,
${\bar\Cal C}'=\Cal C'\otimes_k\bar k$. Then $\Cal T\otimes_k\bar k$ represents

\vskip 5pt
$((\bar k$-schemes$))\to$((sets))

\hskip 52pt $T\mapsto\{$generically smooth $T$-morphisms
$f:T\times_{\bar k}{\bar\Cal C}'\to T\times_{\bar k}\bar{\Cal C}\}$

\vskip 5pt
\flushpar Hence $\Cal T\otimes_k\bar k$ is unramified over $\bar k$
by Case 1. This shows that $\Cal T$ is unramified over $k$.

Case 3: general case. We need to show that $\Omega_{\Cal T/S}^1=0$, for this
it is enough to show $(\Omega_{\Cal T/S}^1)_t=0$ for any closed point $t\in\Cal T$.
Let $s=q(t)\in S$, and let $k$ be the residue field at $s$ (i.e.
$\{ s\}$ can be viewed as a morphism $\text{Spec}(k)\to S$).
Denote by $\Cal C_s,\Cal C_s^{\prime}$ and $\Cal T_s$ the fibers of $\Cal C,\Cal C'$
and $\Cal T$ over $s$ respectively. Then $\Cal T_s$ represents

\vskip 5pt
$((k$-schemes$))\to$((sets))

\hskip 52pt $T\mapsto\{$generically smooth $T$-morphisms
$f:T\times_k\Cal C_s^{\prime}\to T\times_k\Cal C_s\}$

\vskip 5pt
\flushpar Hence $\Cal T_s$ is unramified over $k$ by Case 2, i.e.
$\Omega_{\Cal T_s/k}^1=0$. Note that
$\Omega_{\Cal T_s/k}^1\cong\Omega_{\Cal T/S}^1|_{\Cal T_s}$,
we have $\Omega_{\Cal T/S}^1=0$. \ \ \ Q.E.D.

\vskip 10pt
\flushpar{\bf Lemma 3}. Let $\pi :X\to S$ be an unramified separated morphism of
noetherian schemes. Then $\pi$ has at most a finite number of sections
$\zeta :S\to X$.

\vskip 10pt
\flushpar {\sl Proof.} Let $\zeta :S\to X$ be a section of $\pi$. It is easy
to see the following diagram is cartesian:
$$\CD S @>{\zeta}>> X \cr
@VV{\zeta}V @VV{\Delta}V \cr X @>{\beta}>> X\times_SX \endCD\tag6$$
where $\beta =(\zeta\circ\pi ,\text{id}_X)$ (i.e.
$\beta (x)=(\zeta (\pi (x)),x)$). Hence $\zeta$ is a closed
immersion because $\Delta$ is a closed immersion. Let $\Cal J$ be
the ideal sheaf of $\zeta (S)\subset X$, and $\Cal I$ be the ideal
sheaf of $\Delta (X)\subset X\times_SX$. Then $\beta^*\Cal
I\cong\Cal J$ because (6) is cartesian. Since $\pi$ is unramified,
we have $\Cal I=\Cal I^2$, hence $\Cal J=\Cal J^2$. This shows
that $X$ is a disjoint union of $\zeta (S)$ with another closed
subscheme, hence each connected component of $S$ maps
isomorphically to a connected component of $X$ under $\zeta$. Let
$\goth X$ be the set of connected components of $X$ and $\goth S$
be the set of connected components of $S$ (both are finite sets).
Then any $\zeta$ is uniquely determined by an injective map from
$\goth X$ to $\goth S$. Hence $\pi$ has at most a finite number of
sections. \ \ \ Q.E.D.

\vskip 10pt
\flushpar{\bf Theorem 1}. Let $S$ be a noetherian scheme and $\Cal C,\Cal C'$ be curves
over $S$, where $\Cal C\to S$ has genus $g\ge 2$. Then

\vskip 5pt\leftskip=30pt\parindent=-10pt
i) There is a quasi-finite
unramified $S$-scheme $\Cal T$ representing the following functor

\vskip 5pt
$((S$-schemes$))\to$((sets))

\hskip 53pt $T\mapsto\{$generically smooth $T$-morphisms
$f:T\times_S\Cal C'\to T\times_S\Cal C\}$

\vskip 5pt
ii) There are at most a finite number of generically smooth
$S$-morphisms from $\Cal C'$ to $\Cal C$.

iii) In particular, if $S=\text{Spec}(k)$ for a field $k$, then there is a finite
\'etale $k$-scheme $\Cal T$ representing the following functor

\vskip 5pt
$((k$-schemes$))\to$((sets))

\hskip 52pt $T\mapsto\{$generically smooth $T$-morphisms
$f:T\times_k\Cal C'\to T\times_k\Cal C\}$

\vskip 10pt
\flushpar and there are at most a finite
number of finite separable $k$-morphisms from $\Cal C'$ to $\Cal C$.

\vskip 10pt\leftskip=0pt\parindent=20pt
\flushpar {\sl Proof.} i) For any connected component $U\subset S$, the fibers
of $C'$ over $U$ all have same genus. Let $S_1,...,S_n$ be the connected
components of $S$, $\Cal C_i=\Cal C\times_SS_i$,
$\Cal C_i^{\prime}=\Cal C'\times_SS_i$ ($1\le i\le n$). By Proposition 1,
there is a quasi-finite
unramified $S_i$-scheme $\Cal T_i$ representing the following functor

\vskip 5pt
$((S_i$-schemes$))\to$((sets))

\hskip 54pt $T\mapsto\{$generically smooth $T$-morphisms
$f:T\times_{S_i}\Cal C_i^{\prime}\to T\times_{S_i}\Cal C_i\}$

\vskip 5pt
\flushpar ($1\le i\le n$). Clearly $\Cal T=\coprod\limits_{i=1}^n\Cal T_i$
represents

\vskip 5pt
$\{ S$-schemes$\}\to$((sets))

\hskip 48pt $T\mapsto\{$generically smooth $T$-morphisms
$f:T\times_S\Cal C'\to T\times_S\Cal C\}$

\vskip 5pt
ii) Let $q:\Cal T\to S$ be the projection, which is quasi-finite, hence
separated. By i), a generically smooth $S$-morphism $\Cal C'\to\Cal C$
is equivalent to a section $\zeta :S\to\Cal T$ of $q$. By Lemma 3, $q$ has
at most a finite number of sections, hence
there are at most a finite number of generically smooth
$S$-morphisms $\Cal C'\to\Cal C$.

iii) When $S=\text{Spec}(k)$ for a field $k$, a $k$-scheme is quasi-finite unramified
iff it is finite \'etale, and a $k$-morphism
$\Cal C'\to\Cal C$ is generically smooth iff it is finite separable.
Hence the statements hold by i) and ii). \ \ \ Q.E.D.

\vskip 10pt
This gives Theorem 0.4, Theorem 0.5 and Theorem 0.1.

\vskip 10pt
\flushpar{\bf Corollary 1}. Let $C$ be a curve of genus $g\ge 2$ over a field $k$.

\vskip 5pt\leftskip=30pt\parindent=-10pt
i) For any finitely generated field extension $K\supset k$, there are at most a
finite number of $k$-algebra homomorphisms $\phi :K(C)\to K$ such that
$K\supset\text{im}(\phi )$ is a separably generated extension.
In other words, $C$ has at most a finite number of smooth $K$-points over $k$.

ii) For any variety $X$ over $k$, there are at most a finite number of
generically smooth $k$-morphisms from $X$ to $C$.

\vskip 10pt\leftskip=0pt\parindent=20pt
\flushpar {\sl Proof.} i) Note that there is a one to one correspondence

\vskip 5pt
\flushpar $\{$smooth $K$-points of $C$ over $k\}\leftrightarrow\{$separably generated
$k$-extensions $K(C)\hookrightarrow K\}$

\vskip 5pt
\flushpar Let $k'$ be the algebraic closure of $k$ in $K$ (i.e.
$k'=\{ a\in K|a$ is algebraic over $k\}$), then $k'\supset k$ is a finite
extension. Let $C'=C\otimes_kk'$, then $K(C')\cong K(C)[k']$. Let
$m=|Gal(k'/k)|$, then every field extension
$K(C)\hookrightarrow K$ over $k$ induces $m$ field
extensions $K(C')\hookrightarrow K$. Therefore it is enough to show that there
are at most a finite number of separably generated $k'$-extensions
$K(C')\hookrightarrow K$. Thus we may assume $k=k'$, i.e. $k$ is algebraically
closed in $K$.

Let $n=\text{tr.deg}(K/k)$. Then we can take $n$ subfields $L_1,...,L_n\subset K$
containing $k$ with $\text{tr.deg}(L_i/k)=n-1$ ($1\le i\le n$), such that for each
$i$, $L_i$ is algebraically closed in $K$, $L_i\subset K$ is separably generated,
and $\bigcap\limits_{i=1}^nL_i=k$. Let $C_i=C\otimes_kL_i$ ($1\le i\le n$).
Then $g(C_i)=g(C)=g\ge 2$.
Since $\text{tr.deg}(K/L_i)=1$, There is an $L_i$-curve $C_i^{\prime}$ such that
$K(C_i^{\prime})\cong K$.

For any $k$-homomorphism $\phi :K(C)\to K$
such that $K$ is separably generated over $\text{im}(\phi )$,
there is at least one $i$ such that $\text{im}(\phi )\not\subset L_i$
(hence $\text{im}(\phi )\cap L_i=k$), and $K(C)\otimes_kL_i\to K$ is smooth,
hence induces an $L_i$-homomorphism $\phi':K(C_i)\to K$ such that
$K\supset\text{im}(\phi')$ is a separable extension,
which is equivalent to a finite separable
$L_i$-morphism $C_i^{\prime}\to C_i$. By Theorem 1.iii), there are at most a
finite number of finite separable $L_i$-morphisms from $C_i^{\prime}$ to $C_i$.
Note that $\phi'$ uniquely determines $\phi$,
hence there is a monomorphism

\vskip 5pt
\flushpar $\{$smooth $K$-points of $C$ over
$k\}\hookrightarrow\{\bigcup\limits_{i=1}^n\{$finite separable
$L_i$-morphisms $C_i^{\prime}\to C_i\}$

\vskip 5pt
\flushpar This shows $\{$smooth $K$-points of $C$ over $k\}$ is a finite set.

ii) Let $K=K(X)$, then any generically smooth $k$-morphism from $f:X\to C$
gives a smooth $K$-point $f':\text{Spec}(K)\to C$, and $f$ is uniquely determined
by $f'$. This gives a monomorphism

\vskip 5pt
\flushpar $\{$generically smooth $k$-morphisms $X\to C\}\hookrightarrow\{$smooth
$K$-points of $C$ over $k\}$

\vskip 5pt
\flushpar Hence $\#\{$generically smooth $k$-morphisms $X\to C\}<\infty$.
\ \ \ Q.E.D.

\vskip 10pt
This gives Theorem 0.2 and Theorem 0.3.

\vskip 20pt
\flushpar {\bf 2. Some consequences}

\vskip 10pt
In this section we will prove Theorem 0.6 and Theorem 0.7. First we generalize
Theorem 1.1.

\vskip 10pt
\flushpar{\bf Theorem 1}. Let $S$ be a noetherian scheme and $\Cal C_1,...,\Cal C_m$,
$\Cal C_1^{\prime},...,\Cal C_n^{\prime}$ be curves over $S$, such that for each $i$
($1\le i\le m$), $\Cal C_i\to S$ has genus $g_i\ge 2$. Let
$\Cal T_{i,j}$ ($1\le i\le m$, $1\le j\le n$) be the $S$-scheme representing
the following functor

\vskip 5pt
$((S$-schemes$))\to$((sets))

\hskip 53pt $T\mapsto\{$generically smooth $T$-morphisms
$f:T\times_S\Cal C_j^{\prime}\to T\times_S\Cal C_i\}$

\vskip 5pt
\flushpar as in Theorem 1.1, and let
$f_{i,j}:\Cal T_{i,j}\times_S\Cal C_j^{\prime}\to\Cal T_{i,j}\times_S\Cal C_i$
be the universal morphism. Let
$X=\Cal C_1\times_S\cdots\times_S\Cal C_m$ and
$Y=\Cal C_1^{\prime}\times_S\cdots\times_S\Cal C_n^{\prime}$.
Then there is a quasi-finite unramified
$S$-scheme $\Cal T$ representing the following functor

\vskip 5pt
\hskip -10pt $\goth M\goth o\goth r_{\text{gs}}:((S$-schemes$))\to$((sets))

\hskip 78pt $T\mapsto\{$generically smooth $T$-morphisms
$T\times_SY\to T\times_SX\}$

\vskip 5pt
\flushpar Furthermore,
$$\Cal T\cong\coprod_{\lambda\in I}\Cal T_{1,\lambda (1)}\times_S\cdots
\times_S\Cal T_{m,\lambda (m)} \tag 1$$
where $I$ is the set of all injective maps from $\{ 1,2,...,m\}$
to $\{ 1,2,...,n\}$, and the universal morphism over
$\Cal T_{\lambda}=\Cal T_{1,\lambda (1)}\times_S\cdots\times_S
\Cal T_{m,\lambda (m)}$ is
$$f_{\lambda}=(f_{1,\lambda (1)}\times_S\cdots\times_Sf_{m,\lambda (m)})\circ
(\text{id}_{\Cal T_{\lambda}}\times_Sq_{\lambda}):
\Cal T_{\lambda}\times_SY\to\Cal T_{\lambda}\times_SX \tag 2$$
where $q_{\lambda}=\text{pr}_{\lambda (1)}\times_S\cdots\times_S\text{pr}_{\lambda (m)}:
Y\to\Cal C_{\lambda (1)}^{\prime}\times_S\cdots\times_S\Cal C_{\lambda (m)}^{\prime}$
(i.e. $(x_1,...,x_n)\mapsto (x_{\lambda (1)},...,x_{\lambda (m)})$).
In particular, if $S$ is connected,
then for any generically smooth $S$-morphism $f:Y\to X$, there is a
$\lambda\in I$ and generically smooth $S$-morphisms
$\phi_i:\Cal C_{\lambda (i)}^{\prime}\to\Cal C_i$ ($1\le i\le m$) such that
$f((x_1,...,x_n)=(\phi_1(x_{\lambda (1)}),...,\phi_m(x_{\lambda (m)}))$.

\vskip 10pt
\flushpar {\sl Proof.} By Lemma 1.1, there is a locally quasi-projective
$S$-scheme $\Cal T$ representing $\goth M\goth o\goth r_{\text{gs}}$.
Let $f:\Cal T\times_SY\to\Cal T\times_SX$ be the universal $\Cal T$-morphism.
We now show that the projection $q:\Cal T\to S$ is quasi-finite and unramified.

Casr 1: $m=1$. For each $j$ ($1\le j\le n$), denote by
$$Y_j=\Cal C_1^{\prime}\times_S\cdots\times_S\Cal C_{j-1}^{\prime}
\times_S\Cal C_{j+1}^{\prime}\times_S\cdots\times_S\Cal C_n^{\prime}\tag 3$$
and $p_j:Y\to Y_j$ the projection.
For each $j$ ($1\le j\le n$), $f$ is equivalent to the morphism
$\phi_j=p_j\times_Sf:\Cal T\times_SY\to\Cal T\times_SY_j\times_SX$ over
$\Cal T\times_SY_j$. Let $U_j\subset\Cal T$ be the largest
open subscheme over which $\phi_j$ is generically smooth. Then it is easy to
see that $\phi_j$ is generically smooth over $U_j\times_SY_j$.
Hence there is an induced $S$-morphism $q_j:U_j\times_SY_j\to\Cal T_{1,j}$.
Since $\Cal T_{1,j}$ is quasi-finite over $S$, $q_j$ factors through $U_j$.
In other words, $f|_{U_j\times_SY}$ is equal to the composition of
the projection $\text{pr}_j:U_j\times_SY\to U_j\times_S\Cal C_j^{\prime}$
and the pull-back of $f_{1,j}$ via an $S$-morphism $h_j:U_j\to\Cal T_j$.
By the universality of $\Cal T$, we see $h_j$ is an isomorphism.
Furthermore, for any point $t\in\Cal T$, at least one $\phi_j$ is generically
smooth over $t$, hence $\Cal T=\bigcup\limits_{j=1}^nU_j$. The above
argument also shows that the $U_j$s are disjoint to each other.
Therefore we have
$$\Cal T\cong\coprod_{j=1}^n\Cal T_{1,j} \tag 4$$
over $S$, and the universal morphism over $\Cal T_{1,j}$ is
$$f_{1,j}\circ\text{pr}_j:
\Cal T_{1,j}\times_SY\to\Cal T_{1,j}\times_SX \tag 5$$

Case 2: general case. Let $T$ be a connected $S$-scheme. Then a $T$-morphism
$\phi :T\times_SY\to T\times_SX$ is equivalent to $m$ $T$-morphisms
$\phi_i=\text{pr}_i\circ\phi :T\times_SY\to T\times_S\Cal C_i$. If $\phi_i$
is generically smooth, then by Case 1,
there is a unique $j$ such that $\phi_i=\psi_{ij}\circ\text{pr}_j$,
where $\psi_{ij}:T\times_S\Cal C_j^{\prime}\to T\times_S\Cal C_i$
is the pull-back of $f_{i,j}$ via an $S$-morphism
$T\to\Cal T_{i,j}$. Denote $\lambda (i)=j$. It is easy
to see that $\phi$ is generically smooth iff every $\phi_i$
is generically smooth and $\lambda (i)\ne\lambda (i')$ for any $i\ne i'$.
Thus $\lambda\in I$, and $\phi$ is equal to the pull-back of
$f_{\lambda}$ via a unique $S$-morphism $T\to\Cal T_{\lambda}$.

From this we see that $\Cal T$ is isomorphic to a disjoint union of all
$\Cal T_{\lambda}$s, hence is quasi-finite and unramified over $S$. \ \ \ Q.E.D.

\vskip 10pt
In particular, in the case when $S=\text{Spec}(k)$ for a field $k$, we have

\vskip 10pt
\flushpar{\bf Corollary 1}. Let $C_1,...,C_m$,
$C_1^{\prime},...,C_n^{\prime}$ be curves over a field $k$, with
$g(C_i)=g_i\ge 2$ ($1\le i\le m$), $g(C_j^{\prime})=g_j^{\prime}$ ($1\le j\le n$). Let
$X=C_1\times_k\cdots\times_kC_m$,
$Y=C_1^{\prime}\times_k\cdots\times_kC_n^{\prime}$. Then for any generically
smooth $k$-morphism $f:Y\to X$,
there is an injective map $\lambda :\{ 1,2,...,m\}\to\{ 1,2,...,n\}$ and $m$
finite separable morphisms $f_i:C_{\lambda (i)}^{\prime}\to C_i$ over $k$,
such that $f(x_1,...,x_n)=(f_1(x_{\lambda (1)}),...,f_m(x_{\lambda (m)}))$.

\vskip 10pt
\flushpar{\bf Lemma 1}. Let $\pi :X\to S$ be a quasi-finite morphism of
noetherian schemes. Then there is an integer $M$ such that for any point
$s\in S$, the fiber $X_s$ has degree $\le M$ over $s$.

\vskip 10pt
\flushpar {\sl Proof.} Since we are only concerned with the fiber degrees,
we can assume $S$ is reduced.

We use noetherian induction on $S$, when $X=\emptyset$ there is nothing to prove.

Suppose $X\ne\emptyset$. Take a generic point $\xi\in X$ such that $\zeta =\pi (\xi )$
is not a specialization of any $\pi (x)$ ($x\in X$), hence any point of
$\pi^{-1}(\zeta )$ is a generic point of $X$ (because $\pi$ is quasi-finite).
Take an open neighborhood $U'\subset S$ of $\zeta$ such that the generic
points of $\pi^{-1}(U')$ are all in $\pi^{-1}(\zeta )$.
Let $V\subset S$ be the closure of $\{\zeta\}$, with reduced induced scheme
structure.
Since $\pi$ is quasi-finite, we can take an irreducible open neighborhood
$U\subset V\cap U'$ of $\zeta$ such that $\pi^{-1}(U)\to U$ is finite.
Furthermore, noting that $\pi^{-1}(U)\to U$ is generically flat, we can take
$U$ such that $\pi^{-1}(U)\to U$ is flat, hence the fibers of $\pi^{-1}(U)\to U$
all have degree $d=\deg (\pi^{-1}(\zeta )/\zeta )$.

Take an open subset $U_1\subset U'$ such that $U=U_1\cap V$. Note that
$\pi^{-1}(U)=\pi^{-1}(U_1)$. Let $S'=S-U_1$ with reduced induced scheme structure,
and let $X'=X\times_SS'$. By noetherian induction, there is an integer $M$
such that for any $x\in X'$, the fiber degree $\deg (X_s/s)\le M$. Hence for
any $x\in X$, we have $\deg (X_s/s)\le\max (M,d)$. \ \ \ Q.E.D.

\vskip 10pt
For any $g\ge 0$, there is a ``catalog space of curves of genus $g$'', which
is a quasi-projective scheme $S_g$ over $\Bbb Z$ together with a curve
$\Cal C_g$ over $S_g$ such that for any curve $C$ of genus $g$ over any field $k$,
there is at least one $k$-point $\text{Spec}(k)\to S_g$ over which the fiber
of $\Cal C_g$ is isomorphic to $C$. (We have many choices of $S_g$, and we don't use the
moduli space $\Cal M_g$ of curves of genus $g$ because $\Cal M_g$ is not a
fine moduli space, i.e. there is no universal curve over $\Cal M_g$.) For any $g\ge 2$
and $g'\ge 0$, denote by
$S_{g,g'}=S_g\times S_{g'}$. Then over $S=S_{g,g'}$ there are two curves
$\Cal C=\Cal C_g\times S_{g'}$ and $\Cal C'=S_g\times\Cal C_{g'}$, of genera
$g$ and $g'$ respectively. By Theorem 1.1, there is a quasi-finite unramified
$S$-scheme $\Cal T=\Cal T_{g,g'}$ representing

\vskip 5pt
$\{ S$-schemes$\}\to$((sets))

\hskip 48pt $T\mapsto\{$generically smooth $T$-morphisms
$T\times_S\Cal C'\to T\times_S\Cal C\}$

\vskip 5pt
\flushpar By Lemma 1, there is an integer $M$ such that for any $s\in S$,
the fiber $\Cal T_s$ has degree $\le M$ over $s$. For any field $k$ and any
two curves $C,C'$ over $k$ with $g(C)=g$, $g(C')=g'$, there is a $k$-point
$s:\text{Spec}(k)\to S$ such that the fibers $\Cal C_s\cong C$,
$\Cal C_s^{\prime}\cong C'$ over $k$. By $\deg (\Cal T_s/s)\le M$,
we see there are at most $M$ generically smooth (i.e. finite separable)
$k$-morphisms from $C'$ to $C$. Denote by $M(g,g')=M$, we get

\vskip 10pt
\flushpar{\bf Proposition 1}. For any two curves $C,C'$ over any field
$k$ with $g(C)=g$, $g(C')=g'$, we have

\vskip 5pt
$\#\{$finite separable $k$-morphisms $C'\to C\}\le M(g,g')$.

\vskip 10pt
\flushpar{\bf Remark 1}. We can take $M(g,g')$ to be the smallest integer such
that Proposition 1 holds. In this way we define an integer-valued function of
two integer variables $g\ge 2$ and $g'$. By Hurwitz's Theorem, it is easy
to see that $M(g,g')=0$ when $g'<g$. For the bound of $M(g,g')$, we have the
following conjecture.

\vskip 10pt
\flushpar{\bf Conjecture}. There are constants $a,b\in\Bbb R$ such that
$M(g,g')\le ag'+b$ for any $g\ge 2$ and any $g'$.

\vskip 10pt
The following example gives an evidence of the conjecture.

\vskip 10pt
\flushpar{\bf Example 1}. Let $C,C'$ be curves over a field
$k$ with $g(C)=g(C')=g\ge 2$. Then by Hurwitz's theorem,
any finite separable $k$-morphism $C'\to C$ is an isomorphism.
Hence $\#\{$finite separable $k$-morphisms $C'\to C\}\le |Aut(C/k)|$.
It is well-known that $|Aut(C/k)|\le ag+b$ for some constants $a,b$
(this can be shown using Hurwitz's theorem). In other words, the
conjecture holds when $g=g'$.

\vskip 10pt
\flushpar{\bf Corollary 2}. Let $K\supset k$ be a finitely generated field extension,
and $g\ge 2$ be an integer. Then there is an
integer $M(g,K/k)$ such that for any curve $C$ over $k$ with genus $g$,

\vskip 5pt
$\#\{$smooth $K$-points of $C$ over $k\}\le M(g,K/k)$.

\vskip 5pt
\flushpar Therefore for any variety $X$ over $k$, there is an
integer $M(g,X/k)$ such that for any curve $C$ over $k$ with genus $g$,

\vskip 5pt
$\#\{$generically smooth $k$-morphisms $X\to C\}\le M(g,X/k)$.

\vskip 10pt
\flushpar {\sl Proof.} For simplicity we may assume $k$ is algebraically closed.
Let $n\! =\!\text{tr.deg}(K/k)$.
Look at the proof of Corollary 1.1.i), there can be found $n$
subfields $L_1,...,L_n\subset K$
containing $k$ with $\text{tr.deg}(L_i/k)=n-1$ ($1\le i\le n$), such that for each
$i$, $L_i$ is algebraically closed in $K$, $L_i\subset K$ is separably generated,
and $\bigcap\limits_{i=1}^nL_i=k$. For each $i$, there is an $L_i$-curve $C_i^{\prime}$
such that $K(C_i^{\prime})\cong K$. Let $g_i^{\prime}=g(C_i^{\prime})$.
For any $k$-curve $C$ of genus $g$,
a smooth $K$-point of $C$ over $k$ is equivalent to a finite separable
$L_i$-morphism $C_i^{\prime}\to C\otimes_kL_i$ for some $i$. Hence

\vskip 5pt
$\#\{$smooth $K$-points of $C$ over $k\}\le\sum\limits_{i=1}^nM(g,g_i^{\prime})$.

\vskip 5pt
\flushpar We can take $M(g,K/k)=\sum\limits_{i=1}^nM(g,g_i^{\prime})$.

The last statement can be easily deduced by the first one, as in the proof of
Corollary 1.1.ii). \ \ \ Q.E.D.

\vskip 10pt
Proposition 1 and Corollary 2 together give Theorem 0.6.

By Corollary 1 and Proposition 1 we get

\vskip 10pt
\flushpar{\bf Theorem 2}. Let $C_1,...,C_m$,
$C_1^{\prime},...,C_n^{\prime}$ be curves over a field $k$, with
$g(C_i)=g_i\ge 2$ ($1\le i\le m$), $g(C_j^{\prime})=g_j^{\prime}$ ($1\le j\le n$). Let
$X=C_1\times_k\cdots\times_kC_m$,
$Y=C_1^{\prime}\times_k\cdots\times_kC_n^{\prime}$. Then

\vskip 5pt
$\#\{$generically smooth $k$-morphisms $Y\to X\}\le\sum\limits_{\lambda\in I}
\prod\limits_{i=1}^mM(g_i,g_{\lambda (i)}^{\prime})$.

\vskip 5pt
\flushpar where $I$ is the set of all injective maps from $\{ 1,2,...,m\}$
to $\{ 1,2,...,n\}$.

\vskip 10pt
\flushpar{\bf Corollary 3}. Notation as in Theorem 2. For the field $K=K(Y)$ we have

\vskip 5pt
$\#\{$smooth $K$-points of $X$ over $k\}\le\sum\limits_{\lambda\in I}
\prod\limits_{i=1}^mM(g_i,g_{\lambda (i)}^{\prime})$.

\vskip 10pt
\flushpar {\sl Proof.} It is enough to show that for each $i$ ($1\le i\le m$), a smooth
$k$-morphism $\text{Spec}(K)\to C_i$ is equivalent to a generically smooth
$k$-morphism $Y\to C_i$. For each $j$ ($1\le j\le n$), let
$$Y_j=C_1^{\prime}\times_k\cdots\times_kC_{j-1}^{\prime}\times_kC_{j+1}^{\prime}
\times_k\cdots\times_kC_n^{\prime} \tag 6$$
and let $K_j=K(Y_j)$, viewed as a subfield of $K$. For any smooth $k$-morphism
$\phi :\text{Spec}(K)\to C_i$, there is a $j$ such that $\phi$
is equivalent to a generically smooth $K_j$-morphism
$f:C_j^{\prime}\otimes_kK_i\to C_i\otimes_kK_i$. This is then equivalent to
a $k$-morphism $\zeta :\text{Spec}(K_i)\to\Cal T_{i,j}$ (notation in Theorem 1). Note that
the image of $\zeta$ is a $k$-point, because $Y_i$ has geometrically integral
fibers. This shows that $f=f_0\otimes_kK_i$ for a finite separable $k$-morphism
$f_0:C_j^{\prime}\to C_i$, hence is equivalent to $f_0\circ\text{pr}_j:Y\to C_i$.
\ \ \ Q.E.D.

\vskip 10pt
Theorem 1, Corollary 1, Theorem 2 and Corollary 3 together give Theorem 0.7.

\vskip 10pt
\flushpar{\bf Example 2}. Let $C$ be a curve of genus $g\ge 2$ over a field $k$,
and let $X=C\times_kC$. Denote by $\iota :X\to X$ the morphism by exchanging factors
(i.e. $\iota (x,y)=(y,x)$). By Theorem 1 and Example 1, we see that any finite
separable $k$-morphism $X\to X$ is an isomorphism, and is either equal to
$\sigma\times_k\tau$ for some $\sigma ,\tau\in Aut (C/k)$, or equal to
$(\sigma\times_k\tau)\circ\iota$ for some $\sigma ,\tau\in Aut (C/k)$.
Hence $Aut(X/k)\cong (\Bbb Z/2\Bbb Z)\ltimes (Aut(C/k)\times Aut(C/k))$.
Therefore $|Aut(X/k)|=2|Aut(C/k)|^2\le 2M(g,g)^2$.

\vskip 10pt

\vskip 20pt\parindent=0pt
\flushpar {\bf References}

\vskip 10pt
[B] M. Baker: Geometry over $\bar\Bbb Q$ of small height, Part I. MSRI introductory
workshop on rational and integral points on higher-dimensional varieties (2006)

\vskip 5pt
[G] H. Grauert: {\sl Mordell's Vermutung \"uber rationale Punkte auf algebraischen
Kurven und Funktionenk\"orper}, Publ. Math. I.H.E.S. (1965)

\vskip 5pt
[Hid] Haruzo Hida: {\sl $p$-adic Automorphic Forms on Shimura
Varieties}. Springer Monographs in Mathematics. Springer (2004)

\vskip 5pt
[Hir] J.W.P. Hirschfeld: The number of points on a curve, and applications.
Rendiconti di Matematica, Serie V\! I\! I, Vol. 26, Roma (2006), 13-28

\vskip 5pt
[Hr] E. Hrushovski: The Mordell-Lang conjecture for function field. J. AMS Vol. 9,
No. 3 (1996), 667-690

\vskip 5pt
[L1] Ke-Zheng Li: Actions of group schemes (I). Compositio Math. 80, 55-74 (1991)

\vskip 5pt
[L2] K. Li: Automorphism group schemes of finite field extensions.
Max-Planck-Institut f\"ur Mathematik Preprint Series 2000 (28)

\vskip 5pt
[L3] Ke-Zheng Li: Push-out of schemes and some applications. Max-Planck-Institut
f\"ur Mathematik Preprint Series 2000 (29)

\vskip 5pt
[L4] Ke-Zheng Li: Vector fields and automorphism groups. Algebraic Geometry
Colloqium, Japan (2004), 119-126

\vskip 5pt
[LO] K. Li \&\ F. Oort: {\sl Moduli of Supersingular Abelian Varieties},
LNM 1680. Springer (1998)

\vskip 5pt
[M1] D. Mumford - The structure of the moduli spaces of curves and abelian
varieties. In: {\sl Actes, Congr\`es international math.} (1970), tome 1, 457-465.
Paris: Gauthier-Villars (1971).

\vskip 5pt
[M2] D. Mumford \&\ J. Fogarty - {\sl Geometric Invariant Theory}. 2nd ed.,
Springer-Verlag, Berlin-Heidelberg-New York (1982).

\vskip 5pt
[N] J. Noguchi: A higher dimensional analogue of Mordell's conjecture over
function fields. Math. Ann. 258 (1981), 207-212

\vskip 5pt
[S] P. Samuel: Compl\'ements \`a un article de Hans Grauert sur la conjecture
de Mordell, Publ. Math. I.H.E.S. tome 29 (1966), 55-62

\vskip 5pt
[Vol] J.F. Voloch: Diophantine geometry in characteristic p: a survey \linebreak
(http://www.ma.utexas.edu/users/voloch/surveylatex/surveylatex.html)

\vskip 5pt
[Voj] P. Vojta: Mordell's conjecture over function fields. Invent. Math. 98 (1989),
115-138

\vskip 5pt
[X] K. Xu: On the elements of prime power order in $K_2$ of a number field.
To appaer in Acta Arithmetica.

\vskip 20pt\parindent=0pt
Kezheng Li

Department of Mathematics

Capital Normal University

Beijing 100037, China

e-mail: kzli\@\!\, gucas.ac.cn

\bye